Original Article

# Some Inequalities for the Polar Derivative of Some Classes of Polynomials

Nuttapong Arunrat[1] and Keaitsuda Maneeruk Nakprasit[1*]

[1]Department of Mathematics, Faculty of Science, Khon Kaen University,

Amphoe Muang, Khon Kaen, 40002, Thailand

* Corresponding author, Email address: kmaneeruk@hotmail.com

**Abstract**

In this paper, we investigate an upper bound of the polar derivative of a polynomial of degree $n$

$$p(z) = (z - z_m)^{t_m}(z - z_{m-1})^{t_{m-1}} \cdots (z - z_0)^{t_0} \left( a_0 + \sum_{v=\mu}^{n-(t_m+\cdots+t_0)} a_v z^v \right)$$

where zeros $z_0, \ldots, z_m$ are in $\{z : |z| < 1\}$ and the remaining $n - (t_m + \cdots + t_0)$ zeros are outside $\{z : |z| < k\}$ where $k \geq 1$. Furthermore, we give a lower bound of this polynomial where zeros $z_0, \ldots, z_m$ are outside $\{z : |z| \leq k\}$ and the remaining $n - (t_m + \cdots + t_0)$ zeros are in $\{z : |z| \leq k\}$ where $k \leq 1$.

## 1. Introduction

Let $k$ be a positive real number. We denote $\{z : |z| < k\}$ and $\{z : |z| \leq k\}$ by $D(0, k)$ and $\overline{D(0, k)}$, respectively.

Consider a polynomial $p(z)$ of degree $n$. Bernstein (1926) presented the well-known inequality $\max_{|z|=1} |p'(z)| \leq n \max_{|z|=1} |p(z)|$. This result is sharp for a

polynomial $p(z) = az^n$ where $a$ is a nonzero complex number. Bernstein's inequality is sharp for a special class of polynomials. Not only upper bounds of $\max_{|z|=1} |p'(z)|$ have been studied, but also lower bounds.

Lax (1944) proved the conjecture which was posed by Erdős for $p(z)$ having no zero in $D(0,1)$ that $\max_{|z|=1} |p'(z)| \leq \frac{n}{2} \max_{|z|=1} |p(z)|$. This bound was improved by Aziz and Dawood (1988) who proved that $\max_{|z|=1} |p'(z)| \leq \frac{n}{2} \big[\max_{|z|=1} |p(z)| - \min_{|z|=1} |p(z)|\big]$. Furthermore, the equality holds for $p(z) = az^n + b$ with $|b| \geq |a|$.

Turán (1939) proved that $\max_{|z|=1} |p'(z)| \geq \frac{n}{2} \max_{|z|=1} |p(z)|$ where $p(z)$ is a polynomial having all its zeros in $\overline{D(0,1)}$. The bounds of Lax (1944) and Turán (1939) are sharp for a polynomial which has all of its zeros on $\{z : |z| = 1\}$.

The improvement of this lower bound was presented by Aziz and Dawood (1988) that $\max_{|z|=1} |p'(z)| \geq \frac{n}{2} \big[\max_{|z|=1} |p(z)| + \min_{|z|=1} |p(z)|\big]$. This new bound is sharp for a polynomial $p(z) = az^n + b$ with $|b| \leq |a|$.

Govil (1991) studied a polynomial $p(z)$ of degree $n$ which has no zero in $D(0,k), k \geq 1$. He proved that

$$\max_{|z|=1} |p'(z)| \leq \frac{n}{1+k} \big[\max_{|z|=1} |p(z)| - \min_{|z|=k} |p(z)|\big].$$

Moreover, Govil (1991) also studied a polynomial $p(z)$ of degree $n$ having all its zeros in $\overline{D(0,k)}, k \leq 1$, and proved that

$$\max_{|z|=1} |p'(z)| \geq \frac{n}{1+k} \Big[\max_{|z|=1} |p(z)| + \frac{1}{k^{n-1}} \min_{|z|=k} |p(z)|\Big].$$

Both bounds are sharp and equalities hold for $p(z) = (z+k)^n$.

Aziz and Shah (1997) studied a lower bound of a derivative of a polynomial $p(z) = a_0 + \sum_{v=\mu}^{n} a_v z^v$, $1 \leq \mu \leq n$, which has all of its zeros in $D(0,k), k \leq 1$.

**Theorem A** (Aziz & Shah, 1997) If $p(z) = a_0 + \sum_{v=\mu}^{n} a_v z^v$, $1 \leq \mu \leq n$, is a polynomial of degree $n$ having all of its zeros in $D(0, k)$, then for $k \leq 1$

$$\max_{|z|=1} |p'(z)| \geq \frac{n}{1+k^\mu} \left[ \max_{|z|=1} |p(z)| + \frac{1}{k^{n-\mu}} \min_{|z|=k} |p(z)| \right].$$

Equality holds for $p(z) = (z^\mu + k^\mu)^{\frac{n}{\mu}}$, where $n$ is a multiple of $\mu$.

Although their Theorem stated that $p(z)$ has all of its zeros in $D(0, k)$, $k \leq 1$, it is described in the proof that the result still holds when $p(z)$ has a zero on $\{z : |z| = k\}$. Thus, we restate their theorem as follows.

**Theorem** [Restate Theorem A] If $p(z) = a_0 + \sum_{v=\mu}^{n} a_v z^v$, $1 \leq \mu \leq n$, is a polynomial of degree $n$ having all its zeros in $\overline{D(0, k)}$, then for $k \leq 1$

$$\max_{|z|=1} |p'(z)| \geq \frac{n}{1+k^\mu} \left[ \max_{|z|=1} |p(z)| + \frac{1}{k^{n-\mu}} \min_{|z|=k} |p(z)| \right]. \tag{1.1}$$

Somsuwan and Nakprasit (2013) investigated a lower bound of a polynomial of degree $n$ which has a zero outside $\overline{D(0, k)}$ and the remaining zeros in $D(0, k), k \leq 1$. One of their results is as follows.

**Theorem 1** (Somsuwan & Nakprasit, 2013) Let $k \leq 1$ and $p(z) = (z - z_0)^s \left( a_0 + \sum_{v=\mu}^{n-s} a_v z^v \right)$, $1 \leq \mu \leq n - s$, be a polynomial of degree $n$ having zero $z_0$ outside $\overline{D(0, k)}$ and the remaining $n - s$ zeros in $D(0, k)$. Then

$$\max_{|z|=1} |p'(z)| \geq \left[ \frac{A}{(1+|z_0|)^s} - \frac{s}{(1+|z_0|)} \right] \max_{|z|=1} |p(z)| + \left[ \frac{A}{k^{n-s-\mu}(k+|z_0|)^s} \right] \min_{|z|=k} |p(z)|, \tag{1.2}$$

where $A = \frac{|1-|z_0||^s (n-s)}{(1+k^\mu)}$.

Nakprasit and Somsuwan (2017) investigated an upper bound of a polynomial of degree $n$ which has a zero in $D(0,1)$ and the remaining zeros outside $D(0, k), k \geq 1$. One of their results is as follows.

**Theorem 2** (Nakprasit & Somsuwan, 2017) Let $k \geq 1$ and $p(z)$ be a polynomial of degree $n$ in the form

$$p(z) = (z-z_0)^s \left( a_0 + \sum_{v=\mu}^{n-s} a_v z^v \right), 1 \leq \mu \leq n-s, 0 \leq s \leq n-1.$$

If a zero $z_0$ is in $D(0,1)$ and the remaining $n-s$ zeros are outside $D(0,k)$, then

$$\max_{|z|=1} |p'(z)| \leq \left[ \frac{s}{(1-|z_0|)} + \frac{A}{(1-|z_0|)^s} \right] \max_{|z|=1} |p(z)| - \frac{A}{(k+|z_0|)^s} \min_{|z|=k} |p(z)|, \quad (1.3)$$

where $A = \frac{(1+|z_0|)^{s+1}(n-s)}{(1+k^\mu)(1-|z_0|)}$.

The inequality (1.3) is sharp for a polynomial $p(z) = z^s(z+k)^{n-s}$.

The *polar derivative* of a polynomial $p(z)$ of degree $n$ with respect to a complex number α, denoted by $D_\alpha p(z)$, is defined by $D_\alpha p(z) = np(z) + (\alpha - z)p'(z)$. Note that $D_\alpha p(z)$ generalizes the derivative of a polynomial in the sense that $\lim_{\alpha \to \infty}(D_\alpha p(z)/\alpha) = p'(z)$.

The bounds of $D_\alpha p(z)$ have been studied by many researchers. For example, Aziz and Shah (1997, 1998) studied upper bounds of $\max_{|z|=1}|D_\alpha p(z)|$ where $p(z)$ is a polynomial of degree $n$ having no zero in $D(0,k), k \geq 1$ and $\alpha \in \mathbb{C}$ with $|\alpha| \geq 1$. Aziz and Rather (1998), Dewan, Singh, & Mir (2009), and Govil and McTume (2004) studied lower bounds of $\max_{|z|=1}|D_\alpha p(z)|$ where $p(z)$ is a polynomial of degree $n$ having all of its zeros in $\overline{D(0,k)}, k \leq 1$ and $\alpha \in \mathbb{C}$ with $|\alpha| \geq 1$.

In this paper, we investigate an upper bound of the polar derivative of a polynomial of degree $n$

$$p(z) = (z-z_m)^{t_m}(z-z_{m-1})^{t_{m-1}} \cdots (z-z_0)^{t_0} \left( a_0 + \sum_{v=\mu}^{n-(t_m+\cdots+t_0)} a_v z^v \right)$$

where zeros $z_0, \ldots, z_m$ are in $D(0,1)$ and the remaining $n - (t_m + \cdots + t_0)$ zeros are outside $D(0, k)$ where $k \geq 1$. Furthermore, we give a lower bound of this polynomial where zeros $z_0, \ldots, z_m$ are outside $\overline{D(0, k)}$ and the remaining $n - (t_m + \cdots + t_0)$ zeros are in $\overline{D(0, k)}$ where $k \leq 1$. Consequently, our results generalize the inequalities (1.2) and (1.3).

## 2. Upper bound of a polar derivative of polynomials having at least one zero in $D(0, 1)$.

In this section, we investigate an upper bound of $\max_{|z|=1}|D_\alpha p(z)|$ where $\alpha \in \mathbb{C}$ with $|\alpha| \geq 1$ and $p(z)$ is a polynomial of degree $n$ which has some zeros in $D(0,1)$ and the remaining zeros outside $D(0, k)$ where $k \geq 1$.

For a polynomial $p(z)$ of degree $n$, we define $q(z) = z^n \overline{p(1/\bar{z})}$. Let $z$ be a complex number with $|z| = 1$. It follows from the result of Govil and Rahman (1969) that

$$|p'(z)| + |q'(z)| \leq n \cdot \max_{|z|=1}|p(z)|. \tag{2.1}$$

Furthermore, one can show that

$$|np(z) - zp'(z)| = |q'(z)|. \tag{2.2}$$

**Theorem 3** Let $p(z)$ be a polynomial of degree $n$ in the form

$$p(z) = (z - z_0)^s \left( a_0 + \sum_{v=\mu}^{n-s} a_v z^v \right), 1 \leq \mu \leq n - s, 0 \leq s \leq n - 1.$$

Let $k \geq 1$ and $\alpha \in \mathbb{C}$ with $|\alpha| \geq 1$. If a zero $z_0$ is in $D(0,1)$ and the remaining $n - s$ zeros are outside $D(0, k)$, then

$$\max_{|z|=1}|D_\alpha p(z)| \leq \left[ n + (|\alpha| - 1)\left( \frac{s}{(1-|z_0|)} + \frac{A}{(1-|z_0|)^s} \right) \right] \max_{|z|=1}|p(z)|$$

$$-\left[\frac{(|\alpha|-1)A}{(k+|z_0|)^s}\right]\min_{|z|=k}|p(z)|, \qquad (2.3)$$

where $A = \frac{(1+|z_0|)^{s+1}(n-s)}{(1+k^\mu)(1-|z_0|)}$.

**Proof:** Observe that

$$|D_\alpha p(z)| = |np(z) - (\alpha - z)p'(z)| \le |np(z) - zp'(z)| + |\alpha||p'(z)| \qquad (2.4)$$

Set $q(z) = z^n \overline{p(1/\bar{z})}$. From the relation (2.2), we have $|np(z) - zp'(z)| = |q'(z)|$ for $|z| = 1$.

By substituting this result into (2.4), we obtain for $|z| = 1$ that

$$|D_\alpha p(z)| \le |q'(z)| + |\alpha||p'(z)| = |q'(z)| + |p'(z)| + (|\alpha| - 1)|p'(z)|.$$

Consequently, $\max_{|z|=1}|D_\alpha p(z)| \le n \cdot \max_{|z|=1}|p(z)| + (|\alpha| - 1)\max_{|z|=1}|p'(z)|$

from the relation (2.1).

Theorem 2 implies that

$$\max_{|z|=1}|p'(z)| \le \left[\frac{s}{(1-|z_0|)} + \frac{A}{(1-|z_0|)^s}\right]\max_{|z|=1}|p(z)| - \frac{A}{(k+|z_0|)^s}\min_{|z|=k}|p(z)|,$$

and therefore

$$\max_{|z|=1}|D_\alpha p(z)| \le n\max_{|z|=1}|p(z)| + (|\alpha| - 1)\left[\left(\frac{s}{(1-|z_0|)} + \frac{A}{(1-|z_0|)^s}\right)\max_{|z|=1}|p(z)|\right.$$

$$\left. - \frac{A}{(k+|z_0|)^s}\min_{|z|=k}|p(z)|\right]$$

$$\le \left[n + (|\alpha| - 1)\left(\frac{s}{(1-|z_0|)} + \frac{A}{(1-|z_0|)^s}\right)\right]\max_{|z|=1}|p(z)|$$

$$-\left[\frac{(|\alpha|-1)A}{(k+|z_0|)^s}\right]\min_{|z|=k}|p(z)|,$$

where $A = \frac{(1+|z_0|)^{s+1}(n-s)}{(1+k^\mu)(1-|z_0|)}$.

**Remark 4** Dividing both sides of the inequality (2.3) by $|\alpha|$ and letting $|\alpha| \to \infty$, we get the inequality (1.3) in Theorem 2.

In case $z_0 = 0$, we obtain the following corollary.

**Corollary 5** Let $p(z)$ be a polynomial of degree $n$ in the form

$$p(z) = z^s\left(a_0 + \sum_{v=\mu}^{n-s} a_v z^v\right), 1 \leq \mu \leq n-s, 0 \leq s \leq n-1.$$

Let $k \geq 1$ and $\alpha \in \mathbb{C}$ with $|\alpha| \geq 1$. If all $n-s$ zeros (except a zero at the origin) are outside $D(0, k)$, then

$$\max_{|z|=1}|D_\alpha p(z)| \leq \left[\frac{|\alpha|(n+sk^\mu)+(n-s)k^\mu}{1+k^\mu}\right] \max_{|z|=1}|p(z)|$$

$$- \left[\frac{(|\alpha|-1)(n-s)}{k^s(1+k^\mu)}\right] \min_{|z|=k}|p(z)|. \qquad (2.5)$$

Next, we show that an upper bound in (2.3) is sharp for a polynomial $p(z) = z^s(z+k)^{n-s}$ where $k$ is a real number with $|k| \geq 1$.

One can see that $|D_\alpha p(z)| = |\{z^s[(n-s)k + \alpha n] + \alpha s k z^{s-1}\}(z+k)^{n-s-1}|$.

Note that $(n-s)k + \alpha n > 0$ because $n, k, s \in \mathbb{Z}^+$ and $\alpha \in \mathbb{R}$ with $|\alpha| \geq 1$.

Furthermore, $\max_{|z|=1}|z^s[(n-s)k+\alpha n] + \alpha s k z^{s-1}|$ and $\max_{|z|=1}|(z+k)^{n-s-1}|$ attain at $z = 1$.

These results yield that

$$\max_{|z|=1}|D_\alpha p(z)| = ((n-s)k + \alpha n + \alpha s k)(1+k)^{n-s-1}. \qquad (2.6)$$

The right side of the inequality (2.3) becomes

$$\left[n + (|\alpha|-1)\left(\frac{s}{(1-0)} + \frac{n-s}{(1+k)(1-0)^s}\right)\right] \max_{|z|=1}|p(z)|$$

$$- \left[\frac{(|\alpha|-1)(n-s)}{(1+k)(k+|z_0|)^s}\right] \min_{|z|=k}|p(z)|$$

$$= ((n-s)k + \alpha n + \alpha s k)(1+k)^{n-s-1},$$

which equals $\max_{|z|=1}|D_\alpha p(z)|$ in (2.6).

This means that an upper bound in Theorem 3 is sharp.

Moreover, this polynomial also makes (2.5) an equality, that is, Corollary 5 is sharp.

**Corollary 6** Let $p(z)$ be a polynomial of degree $n$ in the form

$$p(z) = (z - z_1)^{t_1}(z - z_0)^{t_0}\left(a_0 + \sum_{v=\mu}^{n-(t_1+t_0)} a_v z^v\right), 1 \leq \mu \leq n - (t_1 + t_0),$$

$0 \leq t_1 + t_0 \leq n - 1$.

Let $k \geq 1$ and $\alpha \in \mathbb{C}$ with $|\alpha| \geq 1$. If zeros $z_0$ and $z_1$ are in $D(0,1)$ and the remaining $n - (t_1 + t_0)$ zeros are outside $D(0,k)$, then

$$\max_{|z|=1}|D_\alpha p(z)| \leq \left[\frac{t_1(|\alpha|+|z_1|)(1+|z_1|)^{t_1-1}}{(1-|z_1|)^{t_1}} + \frac{(n-t_1)(1+|z_1|)^{t_1}}{(1-|z_1|)^{t_1}}\right.$$

$$\left. + \frac{(|\alpha|-1)(1+|z_1|)^{t_1}}{(1-|z_1|)^{t_1}}\left(\frac{t_0}{(1-|z_0|)} + \frac{A}{(1-|z_0|)^{t_0}}\right)\right] \max_{|z|=1}|p(z)|$$

$$- \left[\frac{(|\alpha|-1)(1+|z_1|)^{t_1}A}{(k+|z_0|)^{t_0}(k+|z_1|)^{t_1}}\right] \min_{|z|=k}|p(z)|,$$

where $A = \frac{(1+|z_0|)^{t_0+1}(n-(t_0+t_1))}{(1+k^\mu)(1-|z_0|)}$.

**Proof:** Let $p(z) = (z - z_1)^{t_1} p_0(z)$ where $p_0(z) = (z - z_0)^{t_0}\left(a_0 + \sum_{v=\mu}^{n-(t_1+t_0)} a_v z^v\right)$ and $\alpha \in \mathbb{C}$ with $|\alpha| \geq 1$.

Then $D_\alpha p(z) = (z - z_1)^{t_1}[D_\alpha p_0(z)] + t_1(\alpha - z_1)(z - z_1)^{t_1-1} p_0(z)$,

and $|D_\alpha p(z)| \leq |z - z_1|^{t_1}|D_\alpha p_0(z)| + t_1|\alpha - z_1||z - z_1|^{t_1-1}|p_0(z)|$.

Since $|z - z_1| \leq |z| + |z_1| = 1 + |z_1|$ and $|\alpha - z_1| \leq |\alpha| + |z_1|$ for $|z| = 1$, we get

$\max_{|z|=1}|D_\alpha p(z)| \leq (1 + |z_1|)^{t_1} \max_{|z|=1}|D_\alpha p_0(z)|$

$$+ t_1(|\alpha| + |z_1|)(1 + |z_1|)^{t_1-1} \max_{|z|=1}|p_0(z)|.$$

Theorem 3 yields that

$$\max_{|z|=1}|D_\alpha p_0(z)| \leq \left[(n - t_1) + (|\alpha| - 1)\left(\frac{t_0}{(1-|z_0|)} + \frac{A}{(1-|z_0|)^{t_0}}\right)\right] \max_{|z|=1}|p_0(z)|$$

$$- \left[\frac{(|\alpha|-1)A}{(k+|z_0|)^{t_0}}\right] \min_{|z|=k}|p_0(z)|,$$

where $A = \frac{(1+|z_0|)^{t_0+1}(n-(t_0+t_1))}{(1+k^\mu)(1-|z_0|)}$.

9Therefore,

$$\max_{|z|=1}|D_\alpha p(z)| \leq [t_1(|\alpha|+|z_1|)(1+|z_1|)^{t_1-1} + (n-t_1)(1+|z_1|)^{t_1}$$

$$+(|\alpha|-1)(1+|z_1|)^{t_1}\left(\frac{t_0}{(1-|z_0|)}+\frac{A}{(1-|z_0|)^{t_0}}\right)\Big]\max_{|z|=1}|p_0(z)|$$

$$-\left[\frac{(|\alpha|-1)(1+|z_1|)^{t_1}A}{(k+|z_0|)^{t_0}}\right]\min_{|z|=k}|p_0(z)|. \tag{2.7}$$

On $|z|=1$, we have $|p_0(z)| = \frac{1}{|z-z_1|^{t_1}}\cdot|p(z)| \leq \frac{1}{(1-|z_1|)^{t_1}}\cdot|p(z)|$.

Consequently,

$$\max_{|z|=1}|p_0(z)| \leq \frac{1}{(1-|z_1|)^{t_1}}\cdot\max_{|z|=1}|p(z)|. \tag{2.8}$$

On $|z|=k$, we have $|p_0(z)| = \frac{1}{|z-z_1|^{t_1}}\cdot|p(z)| \geq \frac{1}{(k+|z_1|)^{t_1}}\cdot|p(z)|$.

Thus,

$$\min_{|z|=k}|p_0(z)| \geq \frac{1}{(k+|z_1|)^{t_1}}\cdot\min_{|z|=k}|p(z)|. \tag{2.9}$$

By substituting (2.8) and (2.9) in (2.7), we obtain that

$$\max_{|z|=1}|D_\alpha p(z)| \leq \left[\frac{t_1(|\alpha|+|z_1|)(1+|z_1|)^{t_1-1}}{(1-|z_1|)^{t_1}} + \frac{(n-t_1)(1+|z_1|)^{t_1}}{(1-|z_1|)^{t_1}}\right.$$

$$\left.+\frac{(|\alpha|-1)(1+|z_1|)^{t_1}}{(1-|z_1|)^{t_1}}\left(\frac{t_0}{(1-|z_0|)}+\frac{A}{(1-|z_0|)^{t_0}}\right)\right]\max_{|z|=1}|p(z)|$$

$$-\left[\frac{(|\alpha|-1)(1+|z_1|)^{t_1}A}{(k+|z_0|)^{t_0}(k+|z_1|)^{t_1}}\right]\min_{|z|=k}|p(z)|,$$

where $A = \frac{(1+|z_0|)^{t_0+1}(n-(t_0+t_1))}{(1+k^\mu)(1-|z_0|)}$.

**Remark 7** Consider a polynomial of degree $n$

$$p(z) = (z-z_m)^{t_m}(z-z_{m-1})^{t_{m-1}}\cdots(z-z_0)^{t_0}\left(a_0 + \sum_{v=\mu}^{n-(t_m+\cdots+t_0)} a_v z^v\right)$$

where zeros $z_0,\ldots,z_m$ are in $D(0,1)$ and the remaining $n-(t_m+\cdots+t_0)$ zeros are outside $D(0,k)$ where $k \geq 1$.

An upper bound of $\max_{|z|=1}|D_\alpha p(z)|$ where $\alpha \in \mathbb{C}$ with $|\alpha| \geq 1$, can be obtained by applying Theorem 3 as in the proof of Corollary 6.

Let $p_0(z) = (z - z_0)^{t_0}\left(a_0 + \sum_{v=\mu}^{n-(t_m+\cdots+t_0)} a_v z^v\right), p_j(z) = (z - z_j)^{t_j} p_{j-1}(z)$, for $1 \leq j \leq m$, and $\alpha \in \mathbb{C}$ with $|\alpha| \geq 1$. Theorem 3 yields an upper bound of $\max_{|z|=1}|D_\alpha p_0(z)|$. Combining this upper bound together with the facts that

$$\max_{|z|=1}|p_0(z)| \leq \frac{1}{(1-|z_1|)^{t_1}} \cdot \max_{|z|=1}|p_1(z)|$$

and

$$\min_{|z|=k}|p_0(z)| \geq \frac{1}{(k+|z_1|)^{t_1}} \cdot \min_{|z|=k}|p_1(z)|,$$

We can obtain an upper bound of $\max_{|z|=1}|D_\alpha p_1(z)|$ as in Corollary 6.

Consequently, an upper bound of $\max_{|z|=1}|D_\alpha p_j(z)|$ for $2 \leq j \leq m$ can be obtained by a similar process by using an upper bound of $\max_{|z|=1}|D_\alpha p_{j-1}(z)|$ and the facts that

$$\max_{|z|=1}|p_{j-1}(z)| \leq \frac{1}{(1-|z_j|)^{t_j}} \cdot \max_{|z|=1}|p_j(z)|$$

and

$$\min_{|z|=k}|p_{j-1}(z)| \geq \frac{1}{(k+|z_j|)^{t_j}} \cdot \min_{|z|=k}|p_j(z)|,$$

for $2 \leq j \leq m$.

## 3. Lower bound of a polar derivative of polynomials having at least one zero outside $D(0, k)$ where $k \leq 1$.

In this section, we investigate a lower bound of $\max_{|z|=1}|D_\alpha p(z)|$ where $\alpha \in \mathbb{C}$ with $|\alpha| \geq 1$ and $p(z)$ is a polynomial of degree $n$ which has some zeros outside $\overline{D(0,k)}, k \leq 1$ and other zeros in $\overline{D(0,k)}$ where $k \leq 1$.

**Theorem 8** Let $p(z)$ be a polynomial of degree $n$ in the form

$$p(z) = (z - z_0)^s \left( a_0 + \sum_{v=\mu}^{n-s} a_v z^v \right), 1 \leq \mu \leq n - s, 0 \leq s \leq n - 1.$$

Let $k \leq 1$ and $\alpha \in \mathbb{C}$ with $|\alpha| \geq 1$. If a zero $z_0$ is outside $\overline{D(0,k)}$ and the remaining $n - s$ zeros are in $\overline{D(0,k)}$, then

$$\max_{|z|=1} |D_\alpha p(z)| \geq \left[ \frac{(|\alpha|-1)A}{(1+|z_0|)^s} - \left( n + \frac{s(|\alpha|+1)}{(1+|z_0|)} \right) \right] \max_{|z|=1} |p(z)|$$

$$+ \left[ \frac{(|\alpha|-1)A}{k^{n-s-\mu}(k+|z_0|)^s} \right] \min_{|z|=k} |p(z)|, \qquad (3.1)$$

where $A = \frac{|1-|z_0||^s (n-s)}{(1+k^\mu)}$.

**Proof:** By setting $\phi(z) = a_0 + \sum_{v=\mu}^{n-s} a_v z^v$, we can rewrite $p(z) = (z - z_0)^s \phi(z)$.

The derivative of $p(z)$ is $p'(z) = (z - z_0)^s \phi'(z) + \phi(z) s (z - z_0)^{s-1}$ and then

$$D_\alpha p(z) = (\alpha - z)(z - z_0)^s \phi'(z) + [n(z - z_0) + s(\alpha - z)](z - z_0)^{s-1} \phi(z).$$

The triangle inequality implies that

$$|D_\alpha p(z)| + |[n(z - z_0) + s(\alpha - z)](z - z_0)^{s-1} \phi(z)| \geq |(\alpha - z)(z - z_0)^s \phi'(z)|.$$

One can see that

$$\max_{|z|=1} |D_\alpha p(z)| \geq \max_{|z|=1} |(\alpha - z)(z - z_0)^s \phi'(z)|$$

$$- \max_{|z|=1} |[n(z - z_0) + s(\alpha - z)](z - z_0)^{s-1} \phi(z)|.$$

For $|z| = 1$, we obtain that $|z - z_0| \leq |z| + |z_0| = 1 + |z_0|$ and $|\alpha| - 1 = |\alpha| - |z| \leq |\alpha - z| \leq |\alpha| + |z| = |\alpha| + 1$.

Since $|z - z_0|^s \geq (|z| - |z_0|)^s = (1 - |z_0|)^s$ for $k < |z_0| < 1$ and $|z - z_0|^s = |z_0 - z|^s \geq (|z_0| - |z|)^s = (|z_0| - 1)^s$ for $|z_0| > 1$, we obtain that $|z - z_0|^s \geq |1 - |z_0||^s$ for $|z_0| > k$.

Consequently,

$$\max_{|z|=1}|D_\alpha p(z)| \geq (|\alpha|-1)\big|1-|z_0|\big|^s \max_{|z|=1}|\phi'(z)|$$

$$-[n(1+|z_0|) + s(|\alpha|+1)](1+|z_0|)^{s-1} \max_{|z|=1}|\phi(z)|.$$

By applying $\phi(z)$ in the inequality (1.1), we have that

$$\max_{|z|=1}|\phi'(z)| \geq \frac{(n-s)}{(1+k^\mu)}\left[\max_{|z|=1}|\phi(z)| + \frac{1}{k^{n-s-\mu}}\min_{|z|=k}|\phi(z)|\right].$$

This implies that

$$\max_{|z|=1}|D_\alpha p(z)| \geq \left[\frac{(|\alpha|-1)|1-|z_0||^s(n-s)}{(1+k^\mu)}\right.$$

$$\left. -[n(1+|z_0|)+s(|\alpha|+1)](1+|z_0|)^{s-1}\right]\max_{|z|=1}|\phi(z)|$$

$$+\left[\frac{(|\alpha|-1)|1-|z_0||^s(n-s)}{k^{n-s-\mu}(1+k^\mu)}\right]\min_{|z|=k}|\phi(z)|.$$

Observe that $|\phi(z)| = \frac{1}{|z-z_0|^s} \cdot |p(z)| \geq \frac{1}{(1+|z_0|)^s} \cdot |p(z)|$ for $|z| = 1$.

Thus, $\max_{|z|=1}|\phi(z)| \geq \frac{1}{(1+|z_0|)^s} \cdot \max_{|z|=1}|p(z)|$.

On $|z| = k$, we have that $|\phi(z)| = \frac{1}{|z-z_0|^s} \cdot |p(z)| \geq \frac{1}{(k+|z_0|)^s} \cdot |p(z)|$.

This implies that $\min_{|z|=k}|\phi(z)| \geq \frac{1}{(k+|z_0|)^s}\min_{|z|=k}|p(z)|$.

Therefore,

$$\max_{|z|=1}|D_\alpha p(z)| \geq \left[\frac{(|\alpha|-1)|1-|z_0||^s(n-s)}{(1+k^\mu)}\right.$$

$$\left. -[n(1+|z_0|)+s(|\alpha|+1)](1+|z_0|)^{s-1}\right]$$

$$\times\left[\frac{1}{(1+|z_0|)^s}\cdot\max_{|z|=1}|p(z)|\right]$$

$$+\left[\frac{(|\alpha|-1)|1-|z_0||^s(n-s)}{k^{n-s-\mu}(1+k^\mu)}\right]\left[\frac{1}{(k+|z_0|)^s}\min_{|z|=k}|p(z)|\right].$$

$$= \left[\frac{(|\alpha|-1)A}{(1+|z_0|)^s} - \left(n + \frac{s(|\alpha|+1)}{(1+|z_0|)}\right)\right]\max_{|z|=1}|p(z)|$$

$$+ \left[\frac{(|\alpha|-1)A}{k^{n-s-\mu}(k+|z_0|)^s}\right]\min_{|z|=k}|p(z)|,$$

where $A = \frac{|1-|z_0||^s (n-s)}{(1+k^\mu)}$.

**Remark 9** (1) Dividing both sides of the inequality (3.1) by $|\alpha|$ and letting $|\alpha| \to \infty$, we get the inequality (1.2) in Theorem 1.

(2) In case $p(z)$ has at least one zero on $\{z : |z| = k\}$, we obtain that $\min_{|z|=k}|p(z)| = 0$. Then

$$\max_{|z|=1}|D_\alpha p(z)| \geq \left[\frac{(|\alpha|-1)A}{(1+|z_0|)^s} - \left(n + \frac{s(|\alpha|+1)}{(1+|z_0|)}\right)\right] \max_{|z|=1}|p(z)|,$$

where $A = \frac{|1-|z_0||^s (n-s)}{(1+k^\mu)}$.

**Corollary 10** Let $p(z)$ be a polynomial of degree $n$ in the form

$$p(z) = (z-z_1)^{t_1}(z-z_0)^{t_0}\left(a_0 + \sum_{v=\mu}^{n-(t_1+t_0)} a_v z^v\right), 1 \leq \mu \leq n - (t_1+t_0),$$

$0 \leq t_1 + t_0 \leq n - 1$.

Let $k \leq 1$ and $\alpha \in \mathbb{C}$ with $|\alpha| \geq 1$. If zeros $z_0$ and $z_1$ are outside $\overline{D(0,k)}$ and the remaining $n - (t_1 + t_0)$ zeros are in $\overline{D(0,k)}$, then

$$\max_{|z|=1}|D_\alpha p(z)| \geq \left[\frac{(|\alpha|-1)|1-|z_1||^{t_1}A}{(1+|z_0|)^{t_0}(1+|z_1|)^{t_1}} - \left(\frac{(n-t_1)|1-|z_1||^{t_1}}{(1+|z_1|)^{t_1}} + \frac{t_0(|\alpha|+1)|1-|z_1||^{t_1}}{(1+|z_0|)(1+|z_1|)^{t_1}}\right)\right.$$

$$\left. - \frac{t_1(|\alpha|+|z_1|)}{(1+|z_1|)}\right]\max_{|z|=1}|p(z)|$$

$$+ \left[\frac{(|\alpha|-1)|1-|z_1||^{t_1}A}{k^{n-(t_1+t_0)-\mu}(k+|z_0|)^{t_0}(k+|z_1|)^{t_1}}\right]\min_{|z|=k}|p(z)|,$$

where $A = \frac{|1-|z_0||^{t_0}(n-(t_1+t_0))}{(1+k^\mu)}$.

**Proof:** Let $p(z) = (z-z_1)^{t_1}p_0(z)$ where $p_0(z) = (z-z_0)^{t_0}\left(a_0 + \sum_{v=\mu}^{n-(t_1+t_0)} a_v z^v\right)$ and $\alpha \in \mathbb{C}$ with $|\alpha| \geq 1$.

Then $D_\alpha p(z) = (z-z_1)^{t_1}[D_\alpha p_0(z)] + t_1(\alpha-z_1)(z-z_1)^{t_1-1}p_0(z)$,

and $|D_\alpha p(z)| \geq |z-z_1|^{t_1}|D_\alpha p_0(z)| - t_1|\alpha-z_1||z-z_1|^{t_1-1}|p_0(z)|$.

4Since $|z - z_1|^{t_1} \geq (|z| - |z_1|)^{t_1} = (1 - |z_1|)^{t_1}$ for $k < |z_1| < 1$ and $|z - z_1|^{t_1} = |z_1 - z|^{t_1} \geq (|z_1| - |z|)^{t_1} = (|z_1| - 1)^{t_1}$ for $|z_1| > 1$, we obtain that $|z - z_1|^{t_1} \geq \left|1 - |z_1|\right|^{t_1}$ for $|z_1| > k$.

For $|z| = 1$, we get that $|z - z_1| \leq |z| + |z_1| = 1 + |z_1|$ and $|\alpha - z_1| \leq |\alpha| + |z_1|$.

By combining these two results, we have that

$$\max_{|z|=1}|D_\alpha p(z)| \geq \left|1 - |z_1|\right|^{t_1} \max_{|z|=1}|D_\alpha p_0(z)|$$
$$- t_1(|\alpha| + |z_1|)(1 + |z_1|)^{t_1 - 1} \max_{|z|=1}|p_0(z)|. \quad (3.2)$$

By applying $p_0(z)$ in Theorem 8, we obtain that

$$\max_{|z|=1}|D_\alpha p_0(z)| \geq \left[\frac{(|\alpha|-1)A}{(1+|z_0|)^{t_0}} - (n - t_1) - \frac{t_0(|\alpha|+1)}{(1+|z_0|)}\right] \max_{|z|=1}|p_0(z)|$$
$$+ \left[\frac{(|\alpha|-1)A}{k^{n-(t_1+t_0)-\mu}(k+|z_0|)^{t_0}}\right] \min_{|z|=k}|p_0(z)|,$$

where $A = \frac{|1-|z_0||^{t_0}(n-(t_1+t_0))}{(1+k^\mu)}$.

By substituting this result into (3.2), we have that

$$\max_{|z|=1}|D_\alpha p(z)| \geq \left[\frac{(|\alpha|-1)|1-|z_1||^{t_1}A}{(1+|z_0|)^{t_0}} - (n - t_1)|1 - |z_1||^{t_1}\right.$$
$$\left. - \frac{t_0(|\alpha|+1)|1-|z_1||^{t_1}}{(1+|z_0|)} - t_1(|\alpha| + |z_1|)(1 + |z_1|)^{t_1 - 1}\right] \max_{|z|=1}|p_0(z)|$$
$$+ \left[\frac{(|\alpha|-1)|1-|z_1||^{t_1}A}{k^{n-(t_1+t_0)-\mu}(k+|z_0|)^{t_0}}\right] \min_{|z|=k}|p_0(z)|.$$

Since $|p_0(z)| = \frac{1}{|z-z_1|^{t_1}} \cdot |p(z)| \geq \frac{1}{(1+|z_1|)^{t_1}} \cdot |p(z)|$ for $|z| = 1$,

$$\max_{|z|=1}|p_0(z)| \geq \frac{1}{(1+|z_1|)^{t_1}} \cdot \max_{|z|=1}|p(z)|.$$

For $|z| = k$, we have $|p_0(z)| = \frac{1}{|z-z_1|^{t_1}} \cdot |p(z)| \geq \frac{1}{(k+|z_1|)^{t_1}} \cdot |p(z)|$ and then

$$\min_{|z|=k}|p_0(z)| \geq \frac{1}{(k+|z_1|)^{t_1}} \cdot \min_{|z|=k}|p(z)|.$$

Consequently,

$$\max_{|z|=1}|D_\alpha p(z)| \geq \left[\frac{(|\alpha|-1)|1-z_1|^{t_1}A}{(1+|z_0|)^{t_0}(1+|z_1|)^{t_1}} - \left(\frac{(n-t_1)|1-z_1|^{t_1}}{(1+|z_1|)^{t_1}} + \frac{t_0(|\alpha|+1)|1-z_1|^{t_1}}{(1+|z_0|)(1+|z_1|)^{t_1}}\right)\right.$$

$$\left. - \frac{t_1(|\alpha|+|z_1|)}{(1+|z_1|)}\right]\max_{|z|=1}|p(z)|$$

$$+ \left[\frac{(|\alpha|-1)|1-z_1|^{t_1}A}{k^{n-(t_1+t_0)-\mu}(k+|z_0|)^{t_0}(k+|z_1|)^{t_1}}\right]\min_{|z|=k}|p(z)|,$$

where $A = \frac{|1-z_0|^{t_0}(n-(t_1+t_0))}{(1+k^\mu)}$.

**Remark 11** Consider a polynomial of degree $n$

$$p(z) = (z-z_m)^{t_m}(z-z_{m-1})^{t_{m-1}} \cdots (z-z_0)^{t_0}\left(a_0 + \sum_{v=\mu}^{n-(t_m+\cdots+t_0)} a_v z^v\right)$$

where zeros $z_0, \ldots, z_m$ are outside $\overline{D(0,k)}$ and the remaining $n - (t_m + \cdots + t_0)$ zeros are in $\overline{D(0,k)}$ where $k \leq 1$.

A lower bound of $\max_{|z|=1}|D_\alpha p(z)|$, where $\alpha \in \mathbb{C}$ with $|\alpha| \geq 1$, can be obtained by applying Theorem 8 as in the proof of Corollary 10.

Let $p_0(z) = (z-z_0)^{t_0}\left(a_0 + \sum_{v=\mu}^{n-(t_m+\cdots+t_0)} a_v z^v\right)$, $p_j(z) = (z-z_j)^{t_j} p_{j-1}(z)$, for $1 \leq j \leq m$. Theorem 8 yields a lower bound of $\max_{|z|=1}|D_\alpha p_0(z)|$. Combining this lower bound together with the facts that

$$\max_{|z|=1}|p_0(z)| \geq \frac{1}{(1+|z_1|)^{t_1}} \cdot \max_{|z|=1}|p_1(z)|$$

and

$$\min_{|z|=k}|p_0(z)| \geq \frac{1}{(k+|z_1|)^{t_1}} \cdot \min_{|z|=k}|p_1(z)|,$$

We can obtain a lower bound of $\max_{|z|=1}|D_\alpha p_1(z)|$ as in Corollary 10. Consequently, a lower bound of $\max_{|z|=1}|D_\alpha p_j(z)|$ for $2 \leq j \leq m$ can be obtained by a similar process by using a lower bound of $\max_{|z|=1}|D_\alpha p_{j-1}(z)|$ and the facts that

$$\max_{|z|=1}|p_{j-1}(z)| \geq \frac{1}{(1+|z_j|)^{t_j}} \cdot \max_{|z|=1}|p_j(z)|$$

and

$$\min_{|z|=k}|p_{j-1}(z)| \geq \frac{1}{(k+|z_j|)^{t_j}} \cdot \min_{|z|=k}|p_j(z)|,$$

for $2 \leq j \leq m$.

**Acknowledgments**

The authors are financially supported by a scholarship under the Research Fund for Supporting Lecturer to Admit High Potential Student to Study and Research on His Expert Program Year 2016 from the Graduate School, Khon Kaen University, Thailand. (Grant no. 591T106).